\documentclass[12pt]{article}
\usepackage{latexsym}
\usepackage{amsfonts}
\makeatletter

\def\@footnotetext#1{\insert\footins{%
\footnotesize
    \interlinepenalty\interfootnotelinepenalty
    \splittopskip\footnotesep
    \splitmaxdepth \dp\strutbox \floatingpenalty \@MM
    \hsize\columnwidth \@parboxrestore
   \edef\@currentlabel{\csname p@footnote\endcsname\@thefnmark}\@makefntext
    {\rule{\z@}{\footnotesep}\ignorespaces
      #1\strut}}}
\def\abstract{\small\quotation{\hskip-\parindent\sc Abstract.}}

\def\classification{\@ifnextchar [{\@xfootnotenext}%
   {\begingroup\let\protect\noexpand
      \xdef\@thefnmark{}\endgroup
    \@footnotetext}}

\title {}

\begin{document}
\classification {{\it 1991 Mathematics
Subject Classification:} Primary 14E09, 14E25; Secondary
13B10, 13B25.\\ 
$\dagger$) Partially supported by 
  CRCG   Grant 25500/301/01.\\
$\ast$) Partially supported by RGC Fundable Grant
344/024/0002.}  

\begin{center}
{\bf \Large   Embeddings of curves in the plane } 

\bigskip

{\bf   Vladimir Shpilrain}$^{\dagger}$ 

\medskip 

 and 
\medskip

 {\bf Jie-Tai Yu}$^{\ast}$
\smallskip

 Department of Mathematics, University of Hong Kong 

Pokfulam Road, Hong Kong 

e-mail: shpil@hkusua.hku.hk, ~yujt@hkusua.hku.hk

\end{center} 
\medskip 

\begin{abstract}
\noindent   Let $ K[x, y]$ be the polynomial algebra in two variables 
 over a field $K$ of characteristic $0$.  In this paper, we
contribute toward a classification of two-variable polynomials by 
classifying (up to  an automorphism of $ K[x, y]$) 
polynomials of the form $ax^n + by^m + \sum_{im+jn \le
mn} c_{ij} x^i y^j$, $~a, b, c_{ij} \in K$ (i.e.,  polynomials whose 
Newton polygon is either 
a triangle or  a line segment). 
 Our classification has 
several applications to the study of embeddings of algebraic 
curves in the plane.  In particular, 
we show that for any $k \ge 2$, there is an irreducible  
curve  with one place at infinity,  
which has at least $k$ 
inequivalent embeddings in  ${\mathbf C}^2$. Also, upon combining our 
method with a well-known theorem of Zaidenberg and Lin, we show 
that one can decide ``almost" just by inspection whether or not 
a polynomial  fiber ${\{}p(x,y)=0{\}}$ is an irreducible simply 
connected curve.
\end{abstract}

\bigskip

\noindent {\bf 1. Introduction }
\bigskip

 Let $ K[x, y]$ be the polynomial algebra in two variables 
 over a field $K$ of characteristic $0$. Here we contribute toward a 
classification of polynomials from $ K[x, y]$ by proving the 
following 
\medskip

\noindent {\bf Theorem 1.1.} Let $p(x,y)= ax^n + by^m + 
\sum_{im+jn \le mn} c_{ij} x^i y^j$, $~a, b, c_{ij} \in K, ~i,j>0$; 
$a, b$ are not both zero, 
and 
$q(x,y)= Ax^r + By^s + 
\sum_{is+jr \le rs} b_{ij} x^i y^j$, $A, B, b_{ij} \in K, ~i,j>0$; 
 $A, B$ are not both zero. 
Suppose that 
$m$ does not divide $n$, $n$ does not divide $m$, $s$ does not
divide  $r$,$~r$ does not divide  $s$,  and 
max$(m,n) \ne \mbox{max}(r,s)$.  Then  there is no 
 automorphism $\alpha \in Aut(K[x, y])$ that takes $p(x,y)$ 
to $q(x,y)$. 
\medskip

 Polynomials of the form given in Theorem 1.1 can be considered  
``canonical models" for their automorphic images. Note that if,
say,  $n$  divides $m$, then the degree of the polynomial $p(x,y)$ 
in the statement of Theorem 1.1 can be reduced by applying an 
automorphism of the form 
$(x \to x + \mu \cdot y^{m/n}; y \to \lambda y)$ with $\mu, \lambda \in K$. 
However, there is no 
guarantee that the resulting polynomial will have the same form. 
 This shows how subtle the situation is. 

In some special cases
though, we can   handle those polynomials with $m$ divisible by $n$
or vice versa. 
 This is possible, for example, if some fiber of a given 
polynomial admits a one-variable 
polynomial parametrization $x=u(t); ~y=v(t)$:
\medskip 

\noindent {\bf Proposition 1.2.} Suppose the fibers 
${\{}p(x,y)=0{\}}$, ${\{}q(x,y)=0{\}}$ of two polynomials 
$p,q \in \mathbf C[x,y]$, admit  one-variable 
polynomial parametrizations. Then one can effectively find out 
(even without knowing the parametrizations) 
if there is an 
automorphism of $\mathbf C[x,y]$ that takes $p$ to $q$.
\medskip 

 In particular, 
if some fiber of a given 
polynomial is an irreducible simply 
connected curve, then, by  a well-known theorem of 
Zaidenberg and Lin \cite{ZL}, this  fiber admits a one-variable 
polynomial parametrization. More precisely, they prove 
 that (in case $K= {\mathbf C}$) every polynomial like
that has  a canonical model of the form $x^k - y^l$ with $(k,l)=1$. 
Upon combining this with our method, we have the following 
\medskip 

\noindent {\bf Theorem 1.3.} Let $p(x,y) \in {\mathbf C}[x,y]$ be 
a polynomial whose fiber ${\{}p(x,y)=0{\}}$ is an irreducible simply 
connected curve. Then some automorphism  of ${\mathbf C}[x,y]$ takes 
$p(x,y)$ to $x^k - y^l$ with $(k,l)=1$, and: 
\smallskip 

\noindent {\bf (a)} max$(k,l) \le \mbox{deg}(p(x,y))$;
\smallskip 

\noindent {\bf (b)} either $k$ or  $l$ divides $\mbox{deg}(p(x,y))$; 
\smallskip 

\noindent {\bf (c)} the  Newton polygon of $p(x,y)$  is either 
a triangle 
or  a line segment, i.e., $p(x,y)$  is of the form 
$ax^n + by^m + {\displaystyle 
\sum_{im+jn \le mn}} c_{ij} x^i y^j$, $m, n \ge 0$. 
 If $m$ does not divide $n$, $n$ does not divide $m$, and $m, n
\ne 0$, then $m=k$  or $l$, 
and $n=l$ or $k$, respectively. Otherwise, either $p(x,y)$  is 
linear, or 
the ``leading" 
 part $ax^n + by^m +  {\displaystyle \sum_{im+jn = mn}} c_{ij} x^i y^j$ ~is a proper 
power of some other polynomial.
\medskip 

 Thus, in many situations it is possible to rule out polynomials 
without irreducible simply 
connected fibers just by inspection. In any case, by Proposition 1.2,
there is 
an  effective procedure for deciding if a given 
polynomial fiber is irreducible and simply 
connected. 
\medskip 

 The next application of our method concerns embeddings of algebraic 
curves in the plane.
\medskip

\noindent {\bf Theorem 1.4.} For any $k \ge 2$, there is an 
irreducible  algebraic 
curve (with one place at infinity) which has at least $k$ 
inequivalent embeddings in the plane ${\mathbf C}^2$. 
\medskip 

 More formally, this means the following. Suppose we have two 
 polynomial fibers ${\{}p(x,y)=0{\}}$ and ${\{}q(x,y)=0{\}}$. 
We say that these curves are {\it isomorphic} if the algebras of 
residue classes ${\mathbf C}[x,y]/\langle p(x,y) \rangle$ and 
${\mathbf C}[x,y]/\langle q(x,y) \rangle$ are isomorphic. Here $\langle p(x,y) \rangle$ denotes the 
 ideal of ${\mathbf C}[x,y]$ generated by $p(x,y)$. 
 On the other hand, we say that these curves (or, rather, 
embeddings of the same curve in  ${\mathbf C}^2$) 
are {\it equivalent} if there is an 
automorphism of ${\mathbf C}^2$ that takes one of them onto the other.
 
 Now our Theorem 1.4 says that there are arbitrarily (but finitely) 
many isomorphic  algebraic 
curves in ${\mathbf C}^2$, all of which 
belong to different orbits under the action of the group 
$Aut({\mathbf C}^2)$. A particular example of a curve like
that would be $y=x^{p_0} - y^{p_1p_2...p_k}$, where 
$p_0, p_1, ..., p_k$ are distinct primes, $p_0 > p_1 p_2 \cdot  ... \cdot p_k$. 
\smallskip 

 Note that by a result of Abhyankar and Singh \cite{AbSi}, 
an irreducible  
curve with one place  at infinity cannot have infinitely many 
inequivalent embeddings in   ${\mathbf C}^2$. 
\smallskip 

 We also note that the first example of an 
irreducible  algebraic 
curve with one place at infinity which has at least  2 
inequivalent embeddings in   ${\mathbf C}^2$, was recently 
claimed in \cite{AbSa}.
\smallskip 

 To conclude the Introduction, we say a few words about our 
general method.  It is a well-known result of Jung and 
van der  Kulk that every automorphism of $ K[x, y]$ is 
a product of elementary and linear automorphisms. 
The main difficulty 
 in finding a canonical model for a given polynomial is 
to prove that one  can  find a sequence of elementary and linear 
automorphisms that would reduce the degree {\it at every 
step}, until it is further irreducible by any elementary 
automorphism. Then this last polynomial, whose degree is 
irreducible, will be  a canonical model. 

 To  arrange that, we use two principal ideas. First, we mimic   
elementary automorphisms of $ K[x, y]$ by ``elementary 
transformations" of $K[t] \times K[t]$. Second, we use 
Whitehead's idea of ``peak reduction" (see e.g. \cite{LS}) to 
arrange a sequence of elementary transformations of 
$K[t] \times K[t]$ so that the maximum 
 degree would decrease 
at every step. This is described in the next Section 2. 

 While the ``peak reduction" always works for elementary 
transformations of $K[t] \times K[t]$, the first part 
(mimicking elementary automorphisms of $ K[x, y]$ by elementary 
transformations of $K[t] \times K[t]$) is where the 
difficulty is. We managed to do that for polynomials of the 
form given in Theorem 1.1, and also for polynomials 
$p(x,y)$ whose fiber ${\{}p(x,y)=0{\}}$ admits a one-variable 
polynomial parametrization $x=u(t); ~y=v(t)$ (i.e., this 
fiber is a rational curve with one place at infinity). 
The latter is used in proving Proposition 1.2 and Theorem 1.3.
 Those parametrizable fibers actually constitute the most 
tractable class of plane algebraic curves. It seems plausible 
that every curve like that has a unique embedding in   
${\mathbf C}^2$. Below we give a high-school version of this 
conjecture. 
\medskip

\noindent {\bf  Conjecture.}  Suppose $K[p(t), q(t)] = 
K[u(t), v(t)]$ for some (one-variable) polynomials 
$p(t), q(t), u(t), v(t).$ Let $\mbox{deg}(p(t))=k; 
~\mbox{deg}(q(t))=l; ~\mbox{deg}(u(t))=m;\\ 
\mbox{deg}(v(t))=n$,   and 
max$(m,n) > \mbox{max}(k,l)$. 
 Then either  $m$ divides $n$, or  $n$ divides $m$.
\medskip

 So far, this was established only in the case where 
$p(t)=t^k; ~q(t)=t^l; ~(k,l)=1$, in the aforementioned paper by Zaidenberg and Lin \cite{ZL}. This generalizes earlier results of 
Abhyankar and Moh \cite{AbMo} and Suzuki \cite{Suzuki}.\\

 \noindent {\bf 2. Elementary automorphisms and  peak reduction}
\bigskip

 It is a well-known result of Jung and 
van der  Kulk that every automorphism of $ K[x, y]$ is 
a product of elementary and linear automorphisms. We give here 
  a somewhat more precise statement which can be found in 
\cite[Theorem 6.8.5]{PMCohn}:
\medskip 

\noindent {\bf  Proposition 2.1.}
Every automorphism of $K[x,y]$ is 
a product of linear automorphisms and automorphisms of the form 
$x \to x + f(y); ~y \to y$. More precisely, if $(g_1,
g_2)$ is an automorphism of   $K[x,y]$  such 
  that  $~\mbox{deg}(g_1) \ge \mbox{deg}(g_2)$, say, then either $(g_1, g_2)$ is a
linear automorphism, or there exists a unique $~\mu \in 
 K^{\ast}$  and a positive integer $~d~$  such 
  that  $~\mbox{deg}(g_1 - \mu g_2^d) < \mbox{deg}(g_1)$. 
\medskip

 Now we are going to consider the direct product $K[t] \times K[t]$ of two copies of the one-variable polynomial algebra over 
 $K$, and introduce the following elementary 
transformations (ET) that can be applied to elements of this 
algebra: 
\smallskip 

\noindent {\bf (ET1)} $(u, v) \longrightarrow (u+ \mu \cdot 
v^k, v)$ for some $\mu \in K^\ast; ~k \ge 2$. 
\smallskip 

\noindent {\bf (ET2)} $(u, v) \longrightarrow 
(u, v+ \mu \cdot u^k)$. 
\smallskip 

\noindent {\bf (ET3)} a non-degenerate linear transformation 
$(u, v) \longrightarrow (a_1u + a_2v, ~b_1u + b_2v);\\
 a_1, a_2, b_1, b_2 \in K$. 
\smallskip 

  One might notice that some of these transformations  are 
redundant, e.g., (ET1)  is a composition of the other ones. There is 
a reason behind that which will be clear a little later. 

 We shall need the following 
\medskip

\noindent {\bf  Proposition 2.2.} For any pair $(u, v) \in K[t] \times K[t]$,  there is a (perhaps, empty) 
sequence of 
elementary transformations that takes $~(u, v)$ ~to some 
$~(\hat{u}, \hat{v})$   ~such  \\
that:
\smallskip 

\noindent {\bf (i)} the maximum of the degrees 
of polynomials decreases {\it at every step} in this sequence; 
\smallskip 

\noindent {\bf (ii)} the maximum of the degrees 
in $(\hat{u}, \hat{v})$ is irreducible by {\it any} sequence of 
elementary transformations. 
\smallskip 

 Comment to {\bf (i)}: if it happens so that $u$ and  $v$ have 
 the same leading terms, 
 then, perhaps by somewhat abusing the language, 
we say that the transformation $(u, v) \to (u-v, v)$ reduces 
the maximum of the degrees. 
\medskip 

\noindent {\bf Proof.} We shall use the ``peak reduction" method 
to prove this statement. This means the following. If at some 
point of a sequence of ET, the maximum degree goes up (or 
remains unchanged) before eventually going down, then there must 
be a pair of {\it subsequent} ET in our sequence (a ``peak")
such 
  that  one of them increases the maximum degree (or leaves it 
unchanged), and then the other one decreases it. We are going to 
show that such a peak can always be reduced. In other words, {\it 
if the maximum degree can be decreased by a sequence of ET, then 
it can also be decreased by a single ET}. 
To prove   that, we have 
to consider many different cases, but all of them are quite simple.

 Let $(u, v)$ be a pair of polynomials from $K[t] \times K[t]$ 
with, say, $\mbox{deg}(u) \le \mbox{deg}(v)$,  and let $\alpha_1$ 
and $\alpha_2$ be two  subsequent ET applied to $(u, v)$, 
as described in the previous 
paragraph. Consider several cases: 
\smallskip 

\noindent {\bf (1)} $\alpha_1 : (u, v) \longrightarrow (u+ \mu \cdot 
v^k, v)$ for some $\mu \in K^\ast; ~k \ge 2$. 
 
 This $\alpha_1$ strictly
 increases the maximum degree since $\mbox{deg}(u) \le \mbox{deg}(v)$ by the assumption.
 Now we have two possibilities for $\alpha_2$ since 
a linear ET cannot decrease the maximum degree. 
\smallskip 

\noindent {\bf (a)} $\alpha_2 : (u+ \mu \cdot v^k, v)\longrightarrow (u+ \mu \cdot v^k, ~v+\lambda (u+ \mu \cdot v^k)^m)$ for some $\lambda  \in K^\ast; ~m \ge 2$.  
But this obviously {\it increases} the maximum degree, contrary 
to our assumption. 
\smallskip 

\noindent {\bf (b)} $\alpha_2 : (u+ \mu \cdot v^k, v)\longrightarrow
  (u+ \mu \cdot v^k + \lambda \cdot v^m, v)$.  If this $\alpha_2$ 
 decreases the maximum degree, then we should have  
$\mu \cdot v^k = -\lambda \cdot v^m$, in which case $\alpha_2=\alpha_1^{-1}$, and the peak reduction is just 
cancelling out $\alpha_1$ and  $\alpha_2$.
\smallskip 

\noindent {\bf (2)} $\alpha_1 : (u, v) \longrightarrow (u, v+ \mu \cdot u^k)$ for some $\mu \in K^\ast; ~k \ge 2$. 

  If this $\alpha_1$ increases the maximum degree, this can only 
happen when $\mbox{deg}(v+ \mu \cdot u^k) = \mbox{deg}(u^k)$, 
in which case we argue exactly as in the case (1). However, 
since $\mbox{deg}(u) \le \mbox{deg}(v)$, it might happen that 
this $\alpha_1$ does not change the maximum degree. Then we 
consider two possibilities for $\alpha_2$:
\smallskip 

\noindent {\bf (a)}  $\alpha_2 : (u, v+ \mu \cdot u^k) \longrightarrow
  (u, v+ \mu \cdot u^k + \lambda \cdot u^m)$.  If this $\alpha_2$ 
 decreases the maximum degree, then we should have $m\ge k$. 
 If  $m=k$, then $\alpha_1 \alpha_2$  is equal to a single ET.
If  $m>k$, then, in order for $\alpha_2$ to decrease the maximum 
degree, we must have $\mbox{deg}(v)$  divisible by $\mbox{deg}(u)$, in which case $\alpha_2$ alone would decrease the maximum degree of $(u, v)$, i.e., we can get rid of $\alpha_1$. 
\smallskip 

\noindent {\bf (b)} $\alpha_2 : (u, v+ \mu \cdot u^k)\longrightarrow (u+ \lambda (v+ \mu \cdot u^k)^m, ~v+ \mu 
\cdot u^k)$. But this $\alpha_2$ can only change the degree 
of the first polynomial in the pair, and  this is not where 
the maximum degree was. 
\smallskip 

\noindent {\bf (3)} $\alpha_1$ is linear, i.e., 
$\alpha_1 : (u, v) \longrightarrow (a_1u + a_2v, b_1u + b_2v); 
 ~a_1, a_2, b_1, b_2 \in K$.  
Again, we have 
two possibilities  for $\alpha_2$.
\smallskip 

\noindent {\bf (a)}  $\alpha_2 : (a_1u + a_2v, b_1u + b_2v) \longrightarrow 
 (a_1u + a_2v, ~b_1u + b_2v + \mu (a_1u + a_2v)^k)$. 
 If $k=1$, then   $\alpha_2$ is linear, and  therefore does not 
change the maximum degree. If $k>1$, then $\alpha_2$  might 
decrease  the maximum degree, but this can only happen if 
 $a_2=0$, in which case we could   decrease  the maximum degree 
of $(u, v)$ by a single ET of the type 
(ET2).
\smallskip 

\noindent {\bf (b)}  the case where $\alpha_2$ is of the type 
(ET1), is completely similar.
 \smallskip 

 Thus, in any of the considered cases, if there is a ``peak", 
then  we can reduce the number 
of  ET in the sequence. An obvious inductive argument completes 
the proof of Proposition 2.2. $\Box$ \\

 \noindent {\bf 3. Proof of  Theorem 1.1} 
\bigskip

  To prove Theorem 1.1, it is clearly sufficient to prove the following 
\medskip

\noindent {\bf  Proposition 3.1.} Let $p(x,y)= ax^n + by^m + 
\sum_{im+jn \le mn} c_{ij} x^i y^j$, $a, b, c_{ij} \in K$,
$a, b$ are not both zero. Suppose that 
$m$ does not divide $n$, and  $n$ does not divide $m$. Then no 
automorphism $\alpha \in Aut(K[x, y])$ can reduce the 
degree of $p(x,y)$. 
\smallskip 

  First, we need the following 
\medskip

\noindent {\bf  Lemma 3.2.}  Let $p(x,y)$ be a polynomial of the  form  $ax^n + by^m + 
\sum_{im+jn \le mn} c_{ij} x^i y^j$, $c_{ij} \in K$.
Then applying an elementary or linear  
automorphism $\beta$ to $p(x,y)$  gives a polynomial of the same 
form, except,   perhaps, in the case  where $m$ divides $n$ 
 or $n$ divides $m$, say,  $m=kn$,  and $\beta : x \to x+ \mu \cdot 
y^k; ~y \to y$ for some $\mu \in K^\ast$. 
\smallskip 

\noindent {\bf Proof.} The statement is obvious for a linear  
automorphism, so suppose we have an elementary automorphism 
$\beta : x \to x+ \mu \cdot y^k; ~y \to y$ for some 
$\mu \in K^\ast; ~k \ge 2$. Then 
\begin{eqnarray}
\beta(p(x,y))= ax^n + a \mu^ny^{kn}+ 
by^m + \sum_{i=1}^{n-1} b_{ij} x^i y^{k(n-i)} ~+  \nonumber \\
+ \sum_{im+jn \le mn} c_{ij} (x^i+ \mu^iy^{ki+j}+ \sum_{s=1}^{i-1} 
a_{ij} x^s y^{k(i-s)+j}).   
\end{eqnarray} 

 Now we have to consider 3 cases:
\smallskip 

\noindent {\bf (a)}  $kn<m$. We have to show that for every 
monomial $x^i y^j$ in (1), one has $im+jn \le mn$. This  is not 
obvious only for monomials of the form  $x^s y^{k(i-s)+j}.$ 
 Compute: 
 \begin{eqnarray}
sm+ (k(i-s)+j)n = sm+ kn(i - s) + jn.  
\end{eqnarray} 

 To see that the right hand side of (2) does not exceed $mn$, 
note that $sm+ kn(i-s) < sm + m(i-s) = mi$, since  $kn<m$.
Therefore, $sm+ (k(i-s)+j)n < mi + nj \le mn$ by the assumption.  
\smallskip 

\noindent {\bf (b)}  $kn>m$. In this case, the ``leading part" 
of $\beta(p(x,y))$ is $ax^n + a \mu^n y^{kn}$, so we have to show
that for  every monomial $x^i y^j$ in (1), one has $ikn+jn \le
kn^2$.  Again, we only have to show that for monomials of the form  
$x^s y^{k(i-s)+j}$: 
 \begin{eqnarray}
skn + (k(i-s)+j)n = kni + jn.
\end{eqnarray} 
 Since $mi + nj \le mn$, after multiplying both sides by $\frac{kn}{m}$ 
we get $kni + \frac{kn^2}{m}j \le kn^2$. Since $\frac{kn^2}{m} > n$
 (recall that  $kn>m$), this yields $kni + jn \le kn^2$. Comparing 
this to (3) completes the proof in this case. 
\smallskip 

\noindent {\bf (c)}  $kn=m$. The same argument as in the 
previous case works here, unless $\beta : x \to x+ \mu \cdot 
y^k; ~y \to y$ for some $\mu \in K^\ast$, 
which can  cause 
cancellation of the leading $y^m$ and loosing control thereby.
 $\Box$ 
\smallskip 

\noindent {\bf Proof of Proposition 3.1.} By way of contradiction, 
assume there is 
 $\alpha \in Aut(K[x, y])$ that takes $p(x,y)$ to some $q(x,y)$ 
of smaller degree. Put into correspondence to the polynomial 
$p(x,y)$ a pair of its {\it face polynomials} $(p(0,t),p(t,0)) 
\in K[t] \times K[t]$. 

 In  the sequence of elementary (linear)  automorphisms that 
corresponds to the automorphism $\alpha$, there must be 
 an elementary automorphism which decreases the degree of the
current  polynomial. Find the {\it first place} in our  sequence of
elementary (linear) automorphisms where we {\it can} apply an
elementary  automorphism which decreases the degree. Suppose this
automorphism is of the form $\beta : x \to x+ \mu \cdot y^k; ~y \to
y$, and it is applied to a polynomial $\tilde p(x,y)$, which we
assume to have the same form as in the statement of  Proposition 3.1
(by Lemma 3.2, we can indeed make this assumption): $\tilde p(x,y)=
\tilde a x^{\tilde n} +  \tilde b y^{\tilde m} + 
\sum_{i \tilde m + j \tilde n \le \tilde m \tilde n} 
\tilde c_{ij} x^i y^j$. 

 If applying the automorphism $\beta$ decreases the degree of the 
polynomial $\tilde p(x,y)$, then, in particular, $\tilde k \tilde n 
= \tilde m$, and  applying the ET of the form 
$(u, v) \longrightarrow (u+ \nu \cdot v^k, v)$ for some $\nu \in
K^\ast$ to the pair of face polynomials $(\tilde p(0,t), \tilde
p(t,0))$, would decrease the maximum of their degrees. (The
coefficient $\nu$ here can be actually computed as  
$\tilde b \tilde a^{-k}$).

 Now Proposition 2.2 implies that there is an ET that decreases the 
maximum of the degrees of the original face polynomials  $(p(0,t),
p(t,0))$. However, given the  hypotheses of Proposition
3.1, this is readily seen to be impossible.  
 This contradiction completes the proof of Proposition 3.1 and  of 
Theorem 1.1 thereby. $\Box$ \\

 \noindent {\bf 4. Embeddings of curves in the plane } 
\bigskip 

 Before we get to the proof of Proposition 1.2 and Theorem 1.3,
 we make a general 
remark. If a polynomial   fiber ${\{}p(x,y)=0{\}}$ admits a 
one-variable polynomial parametrization $x=u(t); ~y=v(t)$, where 
$u(t),v(t)$ have zero constant terms, 
then, by a result of McKay and   Wang \cite{McW}, the polynomial 
$p^m(x,y)$,     where $m = [{\mathbf C}(t):{\mathbf C}(u(t),v(t))]$,
 equals the resultant $R(x,y)=
\mbox{Res}_t(u(t)-x, v(t)-y)$. Moreover, they prove 
\cite[Theorem 5]{McW} that the 
leading  part of $p^k(x,y)$ 
 is obtained the same way (i.e., as a resultant) from the 
  leading  parts of $u(t)$  and  $v(t)$.  This implies,  
in particular, that the  Newton polygon of $p(x,y)$  is either 
a triangle
or  a line segment, i.e., $p(x,y)$  is of the form 
 \begin{eqnarray}
 ax^n + by^m + \sum_{im+jn \le mn} c_{ij} x^i y^j; ~m, n \ge 0  
\end{eqnarray} 
(see \cite[Corollary 6]{McW}). 
 Furthermore,  from the fact that 
$p(x,y)$  is the minimal polynomial for $u(t)$  and  $v(t)$, 
 it follows that, for example,  
$p(x+y^k,y)$  is the minimal polynomial for $u(t)-v^k(t)$  and  
$v(t)$. This establishes a correspondence between elementary
(linear)  automorphisms of $K[x, y]$ applied to $p(x,y)$, and  
elementary (linear) transformations (ET) of  $K[t] \times K[t]$
applied to $(u(t),v(t))$. Theorem 5 of \cite{McW} 
   implies that in any  sequence of elementary (linear) 
automorphisms of $K[x, y]$ applied to $p(x,y)$, all polynomials have 
the form (4),  and, therefore, by our  Proposition 2.2, the
corresponding sequence of  ET  applied to $(u(t),v(t))$ can be
arranged so that it decreases the maximum of the degrees  in a pair
of polynomials at every step.  
\smallskip 

\noindent {\bf Proof of Proposition 1.2.} 
\medskip 

 First we show that both $p$  and  $q$ have 
``canonical models", i.e., automorphic images whose degrees cannot 
be reduced by any automorphisms of $\mathbf C[x, y]$.  Indeed, by the remark
above, both polynomials are of  the form (4). 
If $m$ divides $n$ or  $n$ 
divides $m$, then we can reduce the degree of the polynomial by 
 applying an elementary  automorphism. This  elementary  automorphism
can be easily found: if, say, $kn=m$,  then we apply the 
automorphism $\beta : x \to x+ \mu \cdot 
y^k; ~y \to y$ (see (1)) with $\mu = \sqrt[n]{-b} $.

 Continuing this way, we arrive at a polynomial of  the form (4),
where $m$ does not divide $n$, $n$ does not divide $m$, and, say,
$m<n$. The other polynomial can be reduced to the same form with 
the degrees $m'$ and $n'$, respectively. Again, we can assume that 
$m'<n'$. Now, if $n' \ne n$, we conclude (by Theorem 1.1) that 
there is no 
 automorphism of $\mathbf C[x, y]$ that takes $p(x,y)$ 
to $q(x,y)$. 

 If $n' = n$, then an automorphism taking $p$ 
to $q$ exists  if and only if a combination of a linear  
automorphism with some automorphism of the form 
${\{}x \to x; ~y \to y+ f(x){\}}$, where  $\mbox{deg}(f) < m/n$,
~can take  $p$ to $q$. 

 To figure out if this is possible, we    have to consider 
coefficients of the polynomial $f(x)$ as indeterminates and 
find out if the corresponding system of polynomial equations 
in those indeterminates (together with indeterminates that 
are the coefficients coming from the linear  
automorphism) has a solution over ${\mathbf C}$. To do that, we 
can apply a well-known algorithm that makes use of Gr\"{o}bner bases (see e.g. \cite{AL}). 

 This latter algorithm is pretty slow in general. However, 
there is one 
special case of Proposition  1.2 where we do not have to apply this 
algorithm at all.  This happens when we want to find out if  a given 
  polynomial is 
{\it coordinate}, i.e.,   is an automorphic image of $x$. 
 In that case, if at some point we get a polynomial of the form (4), where $m$ does not divide $n$ and $n$ does not divide $m$, 
then the polynomial is coordinate if and only if 
max$(m,n)=1$;  no further 
 analysis is needed. $\Box$ 
\smallskip 

\noindent {\bf Proof of Theorem 1.3.}  By a result of 
Zaidenberg and Lin \cite{ZL}, some automorphism  of ${\mathbf C}[x,y]$ takes 
$p(x,y)$ to $x^k - y^l$ with $(k,l)=1$. The polynomial fiber 
 ${\{}x^k - y^l=0{\}}$ admits a 
one-variable polynomial parametrization $x=t^l; ~y=t^k$. 
 Therefore, by the remark in the beginning of this section, 
 the polynomial fiber ${\{}p(x,y)=0{\}}$  
admits a parametrization $x=u(t); ~y=v(t)$, and there is 
a sequence of ET that takes $(u(t),v(t))$ to $(t^l, t^k)$, 
so that the maximum of the degrees 
in a pair of polynomials decreases 
 at every step except, perhaps, several terminal  steps that do not 
change the maximum degree. 
This immediately implies parts (a) and  (c) 
of Theorem 1.3. 

 Then, take the last ET in the sequence that 
 decreases  the maximum 
 degree, i.e., after applying this ET we get a  pair of polynomials 
 whose  degrees are ${\{}k,k{\}}$, or ${\{}k,l{\}}$, or ${\{}l,l{\}}$. 
That means, the preceding 
 pair of polynomials   in the sequence  
either has degrees   ${\{}km,k{\}}$, or ${\{}l,ml{\}}$ for  some 
$m \ge 2$. 
 In any case, either $k$ or  $l$ divides 
the maximum of the degrees 
in the preceding  pair of polynomials, and therefore also divides 
the degree of the corresponding two-variable polynomial. 
An obvious inductive argument completes 
the proof of part (b). $\Box$ 

\medskip 

 Now we get to 
\medskip 

\noindent {\bf  Proof of  Theorem 1.4.} We have to exhibit $k$ 
polynomials $f_1,...,f_k$ from ${\mathbf C}[x,y]$ such 
  that ${\mathbf C}[x,y]/\langle f_i \rangle \cong  
{\mathbf C}[x,y]/\langle f_1 \rangle$ for every $i = 1,..., k$, 
but none of those polynomials can be taken to another by an 
automorphism of ${\mathbf C}[x,y]$ (here the symbol $\cong$ means 
``is isomorphic to"). 

 A particular collection of  such polynomials is as follows 
(it is modeled on the corresponding example in combinatorial  
group  theory \cite{McP}): 
\smallskip 

\noindent $f_1=y-x^{p_0} + y^{p_1p_2...p_k}$, where 
$p_0, p_1, ..., p_k$ are distinct primes, 
$p_0 > p_1 p_2 \cdot  ... \cdot p_k$; 
\smallskip 

\noindent $f_2=y-(x^{p_0} - y^{p_2...p_k})^{p_1}$;
\smallskip 

\noindent $f_3=y-(x^{p_0} - y^{p_3...p_k})^{p_1p_2}$;

...

\noindent $f_k=y-(x^{p_0} - y^{p_k})^{p_1p_2...p_{k-1}}$.
\smallskip 

 We are now going to show that the corresponding algebras of 
residue classes are isomorphic. It will be technically more 
convenient to write those algebras of 
residue classes as ``algebras with relations", 
i.e., for example, instead of 
${\mathbf C}[x,y]/\langle f_1 \rangle$ we shall write 
$\langle x,y \mid f_1=0\rangle$. Following is the chain of 
isomorphism-preserving transformations (similar to Tietze 
transformations in group  theory -- see  \cite{LS}) that 
establishes the isomorphism between $\langle x,y \mid f_1=0\rangle$
 and $\langle x,y \mid f_2=0\rangle$: 
\smallskip 

\noindent $\langle x,y \mid y=x^{p_0} - y^{p_1p_2...p_k}\rangle 
~\cong ~\langle x,y,z \mid y=x^{p_0} - y^{p_1p_2...p_k}; 
~z=y^{p_1}\rangle ~\cong$

\noindent $\langle x,y,z \mid y=x^{p_0} - z^{p_2...p_k}; 
z=y^{p_1}\rangle ~\cong \langle x,z \mid z=
(x^{p_0} - z^{p_2...p_k})^{p_1}\rangle ~\cong$

\noindent $\langle x,y \mid y=(x^{p_0} - y^{p_2...p_k})^{p_1}\rangle $. 
\smallskip 

    In a similar way, one can establish the isomorphism between $\langle x,y \mid f_1=0\rangle$
 and $\langle x,y \mid f_i=0\rangle$ for every $i = 2,..., k$. 

 Applying our Theorem 1.1 shows that none of the 
polynomials $f_i$ can be taken to $f_j, ~j\ne i,$
 by an 
automorphism of ${\mathbf C}[x,y]$ (the  restriction  $p_0 > p_1 p_2 \cdot  ... \cdot p_k$ ensures that the  conditions of Theorem 1.1 are 
      satisfied). 

 Finally,  the fact that  the curve ${\{}f_1=0{\}}$  (and hence 
any curve ${\{}f_i=0{\}}$)  
   is irreducible,  is  obvious since $f_1$ is of the form 
 $u(x) + v(y)$ for non-constant polynomials $u, v$ of relatively 
prime  degrees. 
$\Box$ \\

\noindent {\bf Acknowledgement} 
\medskip 

 We are grateful to V.~Lin for  insightful discussions.

\end{document}